\documentclass[12pt]{article}
\usepackage{amsmath,amssymb,latexsym,graphicx}
\newtheorem{theorem}{Theorem}
\newtheorem{lemma}[theorem]{Lemma}

\newtheorem{remark}[theorem]{Remark}

\newtheorem{definition}[theorem]{Definition}
\newtheorem{conjecture}[theorem]{Conjecture}

\def\pf{{\bf Proof }}

\begin{document}
\title{Holomorphic extension from the unit sphere in $\mathbb C^n$ 
into complex lines passing through a finite  set} 
%\thanks{The work of the second author was supported in part by the NSF Grants DMS 9971674 and
%0002195.} 
\author{Mark~ L.~Agranovsky}
\maketitle
\begin{abstract}
Let $a \in \mathbb C^n.$ Denote $\mathcal{L}_a$ the bunch of complex lines containing $a.$ 
It is easy to construct a real-analytic, and even polynomial, function $f$ on the unit complex 
sphere $\partial B^n$ in $\mathbb C^n$ such that $f$ is the boundary value of no holomorphic function in 
the ball $B^n$, but nevertheless,  for any  complex line $L \in \mathcal{L}_a$ the restriction $f\vert_{L \cap \partial B^n}$ continuously extends in $L \cap \overline {B^n}$ as a function, holomorphic  in $L \cap B^n.$ 
We prove that, however, two bunches of complex lines are sufficient for testing global holomorphic extendibility of real-analytic functions: if $a$ and $b$ are two distinct  points in the closed unit ball and $f \in C^{\omega}(\partial B^n)$  admits the one-dimensional holomorphic extension in any complex line  $L \in \mathcal{L}_a \cup \mathcal{L}_b$ 
then $f$ is  the boundary value of a holomorphic function in the unit ball $B^n.$ 
\end{abstract}

\section{Introduction and the main result}\label{S:Intro}
The problem of testing analyticity by analytic extendibility into complex lines was investigated in many articles.
It was observed in  \cite{AV} that
boundary values of holomorphic functions in the unit ball $\mathbb C^n$ can be characterized by analytic extendibility 
along complex lines. In \cite{St} this result was generalized to general domains $D \subset \mathbb C^n.$ 

It turned out that the family of all complex lines is abundant for the characterization of the 
boundary values of holomorphic functions, and 
can be reduced to narrower those. The extended references can be found in  \cite{GS}. 
There are two natural families of complex lines to be studied: 
lines tangent to a given surface and lines intersecting a given set. 
There is a certain progress in the one-dimensional holomorphic extension problem for the families of the first type (defined by the tangency condition), but we are not mentioning the results in this direction, as in this article we are interested in the families of the second type.

Let us start with some notations and terminology.
First of all, everywhere in the article $n \geq 2.$ 
The notation  $A(\partial \Omega),$ \ $\Omega$ is a bounded domain in $\mathbb C^n$, will be used for the algebra of boundary values of holomorphic functions, i.e., 
all functions $f \in C(\partial \Omega)$ such that there exists a function $F \in C(\overline \Omega),$
holomorphic in $\Omega$ and satisfying the condition $F(z)=f(z), z \in \partial \Omega.$

Given a set $V \subset \mathbb C^n$, 
denote $\mathcal L_V$ the set of complex lines $L$ such that $L \cap V \neq \emptyset.$
We will say that the family $\mathcal {L}_V$ is {\it sufficient} (for testing function from $A(\partial \Omega)$)
if whenever $f \in C(\partial \Omega)$ continuously 
extends, for every $L \in \mathcal {L}_V$, from $L \cap \partial \Omega$ to $L \cap \overline \Omega$ as a function holomorphic in $L \cap \Omega$, one has  $f \in A(\partial \Omega).$

In the introduced terminology, the results of \cite{AV}, \cite{St} just claim that the family $\mathcal {L}_V$ is sufficient when the set $V$ is the whole domain, $V=\Omega.$
In  \cite{AS}, a stronger result was obtained: the family $\mathcal L_V$ is sufficient for any open set  $V$ 
(see also \cite{KM} for generalizations in this direction).

However, counting the parameters shows that the expected minimal dimension of sufficient families should be $2n-2.$
Indeed, the functions $f$ to be characterized are defined on $2n-1$-dimensional real sphere and hence depend
on $2n-1$ real variables. The  condition of the holomorphic extendibility in every fixed complex line  is 1-parametric
(see Section \ref{S:concl}). Therefore, in order to have the total number of the parameters $2n-1$ we need  more
$2n-2$ parameters, which is the minimal possible number of parameters for the sufficient families. 

The minimal $(2n-2)$-parameter family of complex lines, given by the intersection condition, are bunches $\mathcal L_a$ of complex lines intersecting a given one-point set $\{a\}.$ 
However, it is easy to show that this family is insufficient. Indeed,  take $\Omega=B^n$ and $a=0.$ Then
the  function $f(z)=|z_1|^2$ is constant on any circle $L\cap \partial B^n,$ and therefore extends holomorphically
(as a constant) in $L \cap B^n,$  but $f \notin A(\partial B^n)$ because it is real-valued and non-constant.
Notice, that in this example, the function $f$ is real-analytic and, moreover, is a polynomial.

Thus, the above example shows that the minimal possible $(2n-2)$-parametric family of the intersection type - the
bunch of lines through one fixed point is insufficient.
However, it turns out, and this is the main result of this article, that adding one more bunch of complex lines
leads already to a sufficient family for real-analytic functions.

\begin{theorem} \label{T:main}
Let $a, \ b \in \overline {B^n}, \ a \neq b.$ 
Let $f$ be a real-analytic function on the unit sphere $\partial B^n$. Suppose that for every complex line $L,$ containing at least one of the points $a, b$  
there exists a function $F_L \in C(L \cap \overline {B^n}),$ holomorphic in $L \cap B^n$ and such that $f(z)=F_L(z)$ for $z \in L \cap \partial B^n.$ Then $f \in A(\partial B^n).$
\end{theorem}

The proof of Theorem \ref{T:main} rests on the recent result of the author \cite{Agr} on meromorphic extensions from chains of circles in the plane.  The scheme of the proof is as follows. We start with the case $n=2.$
By averaging in rotations in $z_2$-plane we reduce the original problem to a problem for the Fourier coefficients regarded as functions in the unit disc $|z_1|<1.$ The condition of the holomorphic extendibility into complex lines transforms to the condition of holomorphic extendibility of the Fourier coefficients in a family of circles in the plane $z_1.$ Then we use the result of \cite{Agr} and arrive at the needed conclusion in 2-dimensional case. The $n$-dimensional case follows by applying the 2-dimensional result to complex 2-dimensional cross-sections. 

\section{Some lemmas }
We start with $n=2.$

First of all, observe, that the condition of Theorem \ref{T:main} is invariant with respect to the
group $M(B^2)$ of complex Moebius automorphisms of the unit ball $B^2$. Therefore, by applying a suitable automorphisms, we can assume that $a$ and $b$ belong to the complex line $z_2=0,$ without loss of generality.

Thus, from now on $a=(a_1,0), b=(b_1,0).$
For $|a_1|<1$ denote
$$H(a_1,r)=\{|\frac{z_1-a_1}{1-\overline a_1 z_1}|=r\} $$
the hyperbolic circle in the unit disc with the hyperbolic center $a$ and
$$\mathcal H_{a_1}=\{H(a_1,r): 0<r<1\}$$
the family of such circles. For $|a_1|=1$, by $H(a_1,r)$ we understand the horicycle through $a_1$, i.e.,
the circle tangent at $a_1$ from inside to the unit circle.
Denote $$\pi_1: \mathbb C^2 \mapsto \mathbb C$$ the orthogonal projection to the first complex coordinate $z_1.$

The following observation belongs to  J. Globevnik:
\begin{lemma}\label{L:projections} We have
$\pi_1 (\{ L \cap \partial B^2, \ L \in \mathcal{L}_{\{a\}} \})=\mathcal H_{a_1},$ i.e., for $a \in B^2$,
the circles 
$L \cap \partial B^2, \ L \in \mathcal{L}_{\{a\}},$ project orthogonally onto  
the hyperbolic circles with the hyperbolic center $a_1.$ If $a \in \partial B^2,$ 
then instead of hyperbolic circles one obtains horicycles through $a.$
\end{lemma}
\pf Remind that we assume that $a_2=0.$ The conformal automorphism 
$$u_{a_1}(z_1)=\frac{z_1-a_1}{1-\overline {a}_1 z_1}$$
extends to a biholomorphic Moebus transformation $U_a$  of the ball $B^2:$ 
$$U_a(z_1,z_2)=( \frac{z_1-a_1}{1-\overline {a}_1 z_1}, \frac{\sqrt{1-|a|^2}z_2}{1-\overline {a}_1 z_1}).$$
This automorphism preserves complex lines,  moves the point $a=(a_1,0)$ to 0 and,correspondingly,
transforms the family $\mathcal{L}_{\{a\}}$ of complex lines containing $a$ into the family $\mathcal {L}_0$ of complex lines containing 0. 

The case $a=0$ is obvious, the circles $L \cap \partial B^2$, where $L$ runs over the complex lines containing 0,
project onto of the family $\mathcal {H}_0$ of circles in the disc $|z_1| \leq 1$ centered at the origin. 
Then we conclude that the projection of  our family is $u_a^{-1}(\mathcal H_0)=\mathcal H_a.$ 
The case $a \in \partial B^2$ is even simpler. Lemma is proved.

\begin{definition}
We say that a function $F$ in the unit disc $|w|<1$ is
$$F(w)=O((1-|w|^2)^k), \ |w| \to 1,$$ 
where $k \in \mathbb Z,$
if $F(w)=h(w)(1-|w|^2)^k$ where $h(w)$ is continuous on $|w| \leq 1$ and has only isolated zeros of finite order on the boundary circle $|w|=1.$
\end{definition}

Expand $f(z_1,z_2), \ (z_1,z_2) \in \partial B^2,$ in the  Fourier series in the polar coordinates $r, \psi$,
where $z_2=re^{i\psi}, r=|z_2| :$
$$f(z_1,z_2)=\sum\limits_{\nu=-\infty}^{\infty}\alpha_{\nu}(z_1, |z_2|)e^{i\nu\psi}.$$
Since $|z_2|=\sqrt{1-|z_1|^2}$ on the sphere, $\alpha_{\nu}$ depends, in fact, only on $z_1$:
$$\alpha_{\nu}(z_1,|z_2|)=A_{\nu}(z_1).$$
Define 
\begin{equation}\label{E:F}
F_{\nu}(z_1):=\frac{ A_{\nu}(z_1) }{ (1-|z_1|^2)^{\frac{\nu}{2}} }.
\end{equation}
Substituting  $\sqrt{1-|z_1|^2}=|z_2|$ on $\partial B^2$ we have

\begin{equation}\label{E:f(z_1,z_2)}
f(z_1,z_2)=\sum\limits _{\nu=-\infty}^{\infty} F_{\nu}(z_1)z_2^{\nu}.
\end{equation}

The singularity at $z_2=0$ for negative $\nu$ cancels due to vanishing of $F_{\nu}(z_1)$ on $|z_1|=1:$
\begin{lemma}\label{L:Fnu} Let $F_{\nu}$ is defined by (\ref{E:F}). Then
 \begin{enumerate}
\item
$F_{\nu}$ is real-analytic in the open disc $|z_1| \leq 1,$ 
\item
$F_{\nu}(z_1)=O((1-|z_1|^2)^{k}), \ |z_1| \to 1,$ for some $k \geq  0,$
\item
if $\nu <0$ then $k>0.$
\end{enumerate}
\end{lemma}.
\pf 
We have
$$A_{\nu}(z_1)e^{i\nu\psi}=\frac{1}{2\pi}\int\limits_0^{2\pi}f(z_1,e^{i\varphi}z_2)e^{-i\nu\varphi}d\psi.$$
Then $f \in C^{\omega}(\partial B^n)$ implies  $A_{\nu}(z_1)e^{i\nu\psi}$ is real-analytic on 
$\partial B^2.$

By the definition (\ref{E:F})  of $F_{\nu},$ 
\begin{equation}\label{E:alpha}
A_{\nu}(z_1)e^{i\nu\psi}=F_{\nu}(z_1) z_2^{\nu}.
\end{equation}
The left hand side is in $C^{\omega}(\partial B^2)$ because it is obtained by averaging of the real-analytic function $f$. Therefore, the right hand side in (\ref{E:alpha}) is real -analytic.  
If $|z_1|<1$ then $z_2 \neq 0,$ due to $|z_1|^2+|z_2|^2=1.$ When $|z_1|<1$ then $z_2 \neq 0.$ Dividing
both sides in (\ref{E:alpha}) by $z_2^{\nu}$  
one obtains $F_{\nu} \in C^{\omega}(|z_1|<1).$ This proves  statement 1.

Now consider the case $|z_1^0|=1,$  $z_2^0=0$. 
Without loss of generality we can assume $Im z_1^0 \neq 0$ and so $|Re z_1^0| <1.$ Then choose local coordinates $(Re z_1, z_2, \overline {z}_2)$ in a neighborhood of the point
$(z_1^0, z_2^0) \in \partial B^2.$ Then the Taylor series for $f$ near $z_0$ can be written as 
$$f(z)=\sum\limits_{\alpha, \beta \geq 0}  c_{\alpha,\beta}(Re z_1) z_2^{\alpha}\overline{z_2}^{\beta},$$ 
where $c_{\alpha,\beta}(Re z_1)$ are real-analytic in a full neighborhood of $Re z_1^0.$
Substitution $z_2=re^{i\psi}$ yields the expression for the $\nu-th$ term in Fourier series:
\begin{equation}\label{E:beta}
A_{\nu}(z_1)e^{i\psi}=\sum\limits_{\alpha-\beta=\nu}c_{\alpha,\beta}r^{\alpha+\beta}e^{i\nu\psi}
=(\sum\limits_{\beta \geq 0}r^{2\beta}c_{\nu+\beta,\beta}(Re z_1))z_2^{\nu}.
\end{equation}                              
In the last expression we have replaced the indices of the summation $\alpha=\beta+\nu.$

Since $r=|z_2|=\sqrt{1-|z_1|^2}$, by the definition \ref{E:F} of $F_{\nu}$ we obtain:
\begin{equation}\label{E:sum}
F_{\nu}(z_1)=\sum\limits_{\beta \geq 0}(1-|z_1|^2)^{\beta}c_{\nu+\beta,\beta}(Re z_1).
\end{equation}
This gives us the needed representation 
$$F_{\nu}=(1-|z_1|)^{k}h(z_1),$$ 
where $\beta=k$ is the index of the first nonzero term in (\ref{E:sum}). 

On $|z_1|=1$ we have 
$$F_{\nu}(z_1)=c_{\nu+k,k}(Re z_1)$$
and there are only isolated zeros of finite order because the function $c_{\nu+k,k}(u)$ is real-analytic
for $u$ near $Re z_1^0.$ Since all the argument is true for arbitrary $z_1^0$ of  modulus 1, the conclusion about zeros on $|z_1|=1$ follows. This proves statement 2.

If $\nu <0$ then $\beta=\alpha - \nu > \alpha \geq 0$ in (\ref{E:beta})  and therefore the order $k$ (the minimal $\beta$ in the sum (\ref{E:sum})) is positive. This is statement 3. 
Lemma \ref{L:Fnu} is proved.

\begin{lemma}\label{L:fourier}
The function $f \in C(\partial B^n)$
admits analytic extension in the complex lines passing through a point $a=(a_1,0) \in B^n$ if and only if
the functions $F_{\nu}(z_1), \ \nu \geq 0,$ defined in Lemma \ref{L:Fnu}, formula (\ref{E:F}), have the property: 
\begin{enumerate}
\item
if $\nu \leq 0$ then $F_{\nu}$ extends inside any circle $H(a,r)$ as an analytic function in the corresponding disc, 
\item
if $\nu >0$ then $F_{\nu}$ extends inside any circle $H(a,r), 0<r<1$ as a
 meromorphic function with the only singular point-a pole at the hyperbolic center $a$ of the order at most $\nu$.
\end{enumerate}
\end{lemma}

\pf
Since $a$ and $b$  belong to the complex line $z_2=0$, the rotated function $f(z_1,e^{i\varphi}z_2)$ possesses
the same property of holomorphic extendibility into complex lines passing through $a$ or through $b$.
Then for each $\nu \in \mathbb Z$
$$\frac{A_{\nu}(z_1)}
{ |z_2|^ {\nu} } z_2^{\nu} = A_{\nu}(z_1)e^{i\nu\psi}=\frac{1}{2\pi}\int\limits_0^{2\pi}f(z_1,e^{i\varphi}z_2)e^{-i\nu\varphi}d\psi$$ 
also has the same extension property. Here $z_2=|z_2|e^{i\psi}.$

It is clear that, vice versa, if all the terms in the Fourier series
have this extension property, then the function $f$ does.

It remains to show that the analytic extendibility in the complex lines
 of the summands in Fourier series is equivalent to the
above formulated properties of $F_{\nu}.$
But since  on the unit sphere we have $|z_2|=\sqrt{1-|z_1|^2},$ then  from the relation (\ref{E:F})

$$\frac{A_{\nu}(z_1)}{|z_2|^{\nu}}z_2^{\nu}=F_{\nu}(z_1)z_2^{\nu}$$ 
and we conclude that
 $F(z_1)z_2^{\nu}$ extends in the complex lines containing $a.$

Let $L$ be such a line, different from $z_1=a$.
Then $L=\{z_2=k(z_1-a_1)\}$ for some complex number $k \neq 0$ and therefore on $\partial B^n$ we have:
$$F_{\nu}(z_1)z_2^{\nu}=F_{\nu}(z_1)k^{\nu}(z_1-a_1)^{\nu}.$$

The function in the right hand side depends only on $z_1$ and hence extends analytically inside
the projection of $L \cap \partial B^n$ which is a hyperbolic circle $H(a,r),$ by  Lemma \ref{L:projections}.
Thus, $f$ extends analytically in $L \cap \partial B^n$ if and only if for all $\nu$ the function
$F_{\nu}(z_1)(z_1-a_1)^ {\nu}$ extends analytically from each circle $H(a,r), a<r<1,$ 
inside the corresponding disc. For $F_{\nu}$ this exactly means the meromorphic extension claimed
in 2. This completes the proof of Lemma.

\section {Proof of Theorem \ref{T:main} in dimension two}\label{S:n=2}

After all the preparations we have done, Theorem \ref{T:main} for the case $n=2$
follows from the result of \cite{Agr}.

In order to formulate this result, as it is stated in \cite{Agr}, we need to introduce the terminology used
in \cite{Agr}. Let $F$ be a function in the unit disc $\Delta.$  We call $F$ {\it regular} if
a) the zero set of $F$ is an analytic set on which $F$ vanishes to a finite order, b) 
$F(w)=O((1-|w|^2)^{\nu}), |w| \to 1,$ for some integer $\nu.$
In our case, the coefficients  $F_{\nu}(z_1)$ are regular by Lemma \ref{L:Fnu}.

Remind that $H(a,r)=\{|\frac{w-a}{1-\overline w}|=1\}, \ |a| <1, $ denotes the hyperbolic circle of radius $r$
with the hyperbolic center $a.$ We will include the case $|a|=1$ and will use the notation $H(a,r)$ for horicycle through $a$, i.e., for the circle of radius $r,$ tangent at $a$ from inside to the unit circle $|w|=1.$

\begin{theorem}(\cite{Agr}, Corollary 4)\label{T:Agr}
Let $F \in C^{\nu}(\overline\Delta)$ be a regular function in the unit disc. 
Let $a, b \in \overline \Delta, \ a \neq b,$ and suppose that $F$ extends from any circle
$H(a,r)$ and $H(b,r), \ 0<r<1,$ as a meromorphic function with the only singular point-a pole, of order at most 
$\nu,$ at the point $a$ or $b$, correspondingly. Then $F$ has the form:
\begin{equation}\label{E:Ag}
F(w)=h_0(w)+ \frac{h_1(w)}{1-|w|^2}+\cdots + \frac{h_{\nu}(w)}{(1-|w|^2)^{\nu}}.
\end{equation}
where $h_j, \ j=0, \cdots, \nu$ are analytic functions in $\Delta.$
\end{theorem}

Lemma \ref{L:Fnu} and Lemma \ref{L:fourier} just say that functions  $F_\nu(z_1)$  satisfy all the conditions of Theorem \ref{T:Agr}. Therefore the function $F:=F_{\nu}$ has the representation of the form
(\ref{E:Ag}).
Moreover, by Lemma \ref{L:Fnu} , $F_{\nu}(z_1)=O((1-|z_1|^2)^k, \ |z_1| \to 1$ with nonnegative $k.$
Therefore, the number $\nu=-k$ in (\ref{E:Ag}) is nonpositive and hence either $\nu=0$ and $F=F_{\nu}=h_0$ is holomorphic,
or $\nu <0$ and $F_{\nu}=0.$
The same can be explained as follows:  $F_{\nu}$ is continuous for $|z_1| \leq 1$ and hence in the decomposition \ref{E:Ag} (with $w=z_1$) nonzero terms with the negative powers of $1-|z_1|^2$ are impossible and only the first term survives:
$$F_{\nu}=h_{0}$$ 
Therefore, $F_{\nu}$ is analytic in $|z_1|<1.$
Moreover, we know from Lemma \ref{L:Fnu}, statement 3, that if $\nu <0$ then $k>0$ and hence $F_{\nu}(z_1)=0$ for $|z_1|=1.$ By the uniqueness theorem this implies $F_{\nu}=0, \nu <0.$ 

Then from (\ref{E:f(z_1,z_2)}) we have
\begin{equation}\label{E:expan}
f(z_1,z_2)=\sum\limits_{\nu \geq 0}F_{\nu}(z_1)z_2^{\nu}, \ |z_1|^2+|z_2|^2=1, 
\end{equation}
and the right hand side is holomorphic in $z_1,z_2.$ This is just the desired global holomorphic
extension of $f$ inside the ball $B^2.$ Theorem \ref{T:main} for the case $n=2$ is proved.

\begin{remark}\label{R:kernel}
Each function of the form (\ref{E:expan}), where $F_{\nu}$ are of the form (\ref{E:Ag}),
possesses holomorphic extension from the sphere  $\partial B^2$ in all cross-sections of the ball $B^2$ by 
complex lines intersecting the complex line $z_2=0.$ 
First example of the function with this property was given by  J. Globevnik: 
$f(z)=\frac{z_2^k}{\overline z_2}$. This function can be reduced on the sphere to the form (\ref{E:Ag}), (\ref{E:expan}):
$f(z)=\frac{z_2^{k+1}}{1-|z|^2}.$ Notice that $f$ is a function of  finite order of smoothness. 
\end{remark}
\section{Proof of Theorem \ref{T:main} in arbitrary dimension}\label{S:arbitrary_n}
Let $Q$ be arbitrary complex two-dimensional plane in $\mathbb C^n,$ containing the points $a$ and $b$.
Then $f\vert_{Q \cap \partial B^n}$ extends holomorphically in any complex line $L \subset Q$ passing through $a$ or $b$ and Theorem \ref{T:main} for the case $n=2$ implies that $f$ extends as a holomorphic function $f_Q$
in the 2-dimensional cross-section $Q \cap B^n.$ 

Let $L_{a,b}$ be the complex line containing $a$ and $b.$ Then for any two 2-planes $Q_1$ and $Q_2,$
containing $a$ and $b,$  we have $Q_1 \cap Q_2 =L_{a,b}$ and $f_{Q_1}(z)=f_{Q_2}(z)=f(z)$ for $z$ from the closed curve $L_{a,b} \cap \partial B^n.$  Hence by the uniqueness theorem $f_{Q_1}(z)=f_{Q_2}(z)$ in $L_{a,b} \cap B^n.$ 

Thus, the holomorphic extension in the 2-planes $Q$ agree on the intersections and therefore define a function $F$ in $\overline {B^n}.$ This function has the following properties: it is holomorphic on 2-planes $Q, \ a,b \in Q,$ and   $F\vert_{\partial B^n}=f.$ It follows from the construction that $F$ is $C^{\infty}$ (and even real-analytic).

From the very beginning, applying  a suitable Moebius automorphism  of $B^n$, 
we can assume that both  points, $a$ and $b$ , belong to a coordinate complex line, so that $0 \in L_{a,b}.$  
Then $F$ is holomorphic on any complex line
passing through $0.$ By Forelli theorem \cite{R},4.4.5, $F$ is holomorphic in $B^n$ (for real-analytic $F$ holomorphicity of $F$ can be proved directly). The required global holomorphic extension is built.

\section{Concluding remarks}\label{S:concl}
\begin{itemize}
\item
The examples in Remark \ref{R:kernel} show that without strong smoothness assumptions (in our case, real-analyticity)
even the complex lines meeting the set of $n$ affinely independent points are not enough for the one-dimensional holomorphic extension property.  It is natural to conjecture that $n+1$ affinely independent points are enough and not only for the case of complex ball:
\begin{conjecture} Let $\Omega$ be a domain in $\mathbb C^n$ with a smooth boundary and  \\
$a_1,\cdots, a_{n+1} \in \overline \Omega$ are $n+1$ points  belonging to no complex hyperplane.
Suppose that $f \in C(\partial \Omega)$ extends holomorphically in each cross-section $ L \cap \Omega.$
Then $f \in A(\partial \Omega).$ 
\end{conjecture}

To confirm this conjecture for the complex ball $B^n$, it would be enough to prove the following
\begin{conjecture}
Let $a,b  \in \overline B^2, \ a \neq b.$ Any function 
$f \in C(\partial B^2)$ having one-dimensional holomorphic extension into complex lines passing through $a$ or $b$ has the form  (\ref{E:expan}),(\ref{E:Ag}):
$$f(z)=\sum\limits_{\nu=0}^{\infty}F_{\nu}(z_1)z_2^{\nu}, 
\ F_{\nu}=\sum\limits_{j=0}^{\nu}\frac{h_j(z_1)}{(1-|z_1|^2)^{j}},$$
where $h_j(z_1)$ are holomorphic.
\end{conjecture}
%\begin {remark}\label{R:Morera}
%\begin{enumerate}
\item
Theorem \ref{T:main} can be regarded as a boundary Morera theorem (cf. \cite{GS}), because the condition of the holomorphic extendibility into complex lines $L$ can be written as the complex moment condition 
$$\int_{L \cap \partial B^n}f \omega=0,$$
for any holomorphic differential (1,0)-form $\omega.$
For each fixed complex line $L$, the above condition depends on one real parameter, because is suffices to
take $\omega$ such that $\omega\vert_{L}=\frac{d\zeta}{\zeta-t}$, 
where $\zeta$ is the complex parameter on $L$ and $t$ runs over any fixed real curve in $L \setminus \overline {B^n}.$

Results of the type of Theorem \ref{T:main} can be also viewed as boundary analogs of Hartogs' theorem about separate analyticity
(see \cite{AS}).
\item 
In Theorem \ref{T:main}, the points $a, b$ are taken from the closed ball $\overline{B^n}.$ 
It is shown in \cite{KM1} that finite sets $V \subset \mathbb C^n \setminus \overline{B^n}$ can produce insufficient families $\mathcal{L}_V$ of complex lines. 
For example, it is easy to see that the function $f(z)=|z_1|^2$ extends holomorphically (as a constant) from the unit complex sphere in each complex line parallel to a coordinate line.  This family of lines can be viewed as $\mathcal{L}_V$ where $V$ is $n$ points in $\mathbb {CP}^n \setminus \mathbb C^n$. To obtain finite points, one can  apply Moebius transformation since they preserve complex lines. 
%\item
%In Section \ref{S:arbitrary_n} we have derived the $n$-dimensional version of Theorem \ref{T:main}
%from the case $n=2,$ using tangential $CR$ equation.
%A. Tumanov suggested an alernative, more geometric, argument. Briefly, it is as follows.
%The holomorphic extensions in the complex lines passing through $a$ define a function $F_a$ in the punctured ball %$\overline{B^n} \setminus \{a\}.$ Analogously, the holomorphic extensions in the 
%complex lines passing through $b$ define
%the function $F_b$ in $\overline{B^n} \setminus \{b\}.$ It is not hard to prove that the functions $F_a$ and $F_b$ are %smooth (real-analytic) in the corresponding punctured balls. Also  $F_a(z)=F_b(z)$
%for any point  $z \in \overline {B^n} \setminus \{a,b\}$ because by the two-dimensional version of Theorem %\ref{T:main} 
%the function  $f$ admits holomorphic extension into the complex 2-plane $Q$
%containing the three points $a,b,z$ and hence $F_a(z)$ and $F_b(z)$ are just the value of this 2-dimensional %holomorphic extension. Therefore, $F_a$ is smooth (real-analytic) in $B^n.$ In addition, the function $F_a$ is %holomorphic along any complex line passing through $a$ and hence it is holomorphic in any ball $B(a,r) \subset B^n$, %which follows directly from the power series representation or from Forelli theorem \cite{R}. 
%By the real-analyticty, $F_a$ is holomorphic in the entire ball $B^n.$
%It remains to notice that $F_a\vert_{\partial B^n}=f$ by the construction. 
%The desired holomorphic extension is built.  The author thanks A. Tumanov for presenting this proof.

\end{itemize}
\bigskip

After this article was written,  Josip Globevnik informed the author that he
has proved Theorem \ref{T:main} in the case $n=2$ for infinitely smooth functions. The author thanks J. Globevnik
and A. Tumanov for useful discussions of the results of this article and remarks.

\section*{Acknowledments}
This work was partially supported by the grant from ISF (Israel Science Foundatrion) 688/08.
Some of this research was done as a part of European Science Foundation Networking Program HCAA.


\begin{thebibliography}{15}
\bibitem{Agr} M. L. Agranovsky {\em Characterization of polyanalytic functions by meromorphic extensions into chains of circles}, Preprint, 2009, http://arxiv.org/abs/0910.3578v1
\bibitem{AV} M. L. Agranovsky, R.E. Val'sky, {Maximality of invariant algebras of functions},
Siberian Math. Journal, 12 (1971), 1-7.
\bibitem{A0} M. Agranovsky, {\em $CR$ foliations, the strip-problem and Globevnik-Stout conjecture},
C. R. Acad. Sci. Paris, Ser. I 343 (2006), 91-94.
\bibitem{A1} M.Agranovsky, {\em Propagation of boundary $CR$-foliations and Morera type theorems for manifolds
with attached analytic discs}, Advances in Math., 211, 1, (2007), \ 284-326.
\bibitem{A2} M. Agranovsky, {\em $CR$ foliations, the strip-problem
and Globevnik-Stout conjecture}, C. R. Acad. Sci. Paris, Ser.I 343 (2006),
91-94.
\bibitem{A3} M. Agranovsky, {\em Complex dimensions of real manifolds, attached analytic discs and parametric argument principle}, preprint, arXiv:math 0704.2871
%\bibitem{Agr} M. Agranovsky {\em Characterization of polyanalytic functions by meromorphic extensions
%in chains of circles}, Preprint, 2009.
\bibitem{AG} M. Agranovsky and J. Globevnik, {\em Analyticity on circles for rational and real-analytic
functions of two real variables}, J. D'Analyse Math., 91 (2003), 31-65.
\bibitem{AS} M. Agranovsky and A. M. Semenov {\em Boundary analogues of the Hartogs' theorem},
Siberian Math. J., 32 (1991), 1, 137-139.
Trans. Amer. Math. Soc. 280 (1983), 247-254.
\bibitem{G1} J. Globevnik, {\em Analyticity on rotation invariant families
of circles}, Trans. Amer. Math. Soc., 280 (1983), 247-254.
\bibitem{G2} J.Globevnik, {\em Testing analyticity on rotation invariant
families of curves}, Trans. Amer. Math. Soc. 306 (1988), 401-410.
\bibitem{G3} J. Globevnik, {\em Holomorphic extensions and
rotation invariance}, Complex Variables, 24 (1993), 49-51.
\bibitem{G4} J. Globevnik, {\em Analyticity of functions analytic on circles}, to appear
in Journ. Math. Anal. Appl., arXiv:0906.1295.
\bibitem{GG} J. Globevnik, {\em Analyticity on translates of Jordan curves}, Trans. Amer. Math. Soc. 359 (2007), 5555-5565.
\bibitem{GS} J. Globevnik and E. L. Stout, {\em Boundary Morera theorems for holomorphic functions of several complex variables}, Duke Math.J. 64 (1991), 3, 571-615.
\bibitem{KM} A.M. Kytmanov and S. G. Myslivets, {\em Higher-dimensional boundary analogs of the Morera theorem in problems of analytic continuation of functions}, Complex analysis, J. Math. Sci. (N.Y.), 120 (2004), 1842-1867.
\bibitem{KM1} A.M. Kytmanov and S. G. Myslivets, {\em On families of complex lines, sufficient for holomorphic extension}, Math. Notes, 83,4 (2008), 500-505.
\bibitem{R} W. Rudin, {\em Function theory in the unit ball of $\mathbb C^n$}, Springer Verlag, Berlin, Heidelberg, New York, 1980.
\bibitem{St} E. L. Stout, {\em The boundary values of holomorphic  functions of several complex variables},
Duke Math. J., 44, (1977), 1, 105-108.
\bibitem{T1} A.Tumanov,{\em A Morera type theorem in the strip},
Math.Res. Lett., 11 (2004), no. 1, 23-2
\bibitem{T2} A.Tumanov, {\em Testing analyticity on circles},  Amer. J. of Math. 129,3, (2007),785-790. 
\bibitem{Z1} L. Zalcman,{\em Analyticity and the Pompeiu problem},
Arch. Rat. Mech. Anal. 47 (1972), 237-254.
\end{thebibliography}
\end{document}